\documentclass[11pt]{article}
\usepackage{amssymb}\def\msy{\mathbb}

\def\bbbc{{\msy C}}
\def\bbbq{{\msy Q}}

\def\bbbz{{\msy Z}}

\def\deg{{\rm deg}}
\def\mod#1{\ ({\rm mod}\ #1)}

\newtheorem{theorem}[subsection]{Theorem}

\newtheorem{lemma}[subsection]{Lemma}

\begin{document}
\title{An implementation of Runge's method for Diophantine equations}
\author{F.Beukers and Sz.Tengely}
\maketitle

\begin{abstract} In this paper we suggest an implementation of
Runge's method for solving Diophantine equations satisfying
Runge's condition. In this implementation we avoid the use of
Puiseux series and algebraic coefficients.
\end{abstract}

\section{Introduction}
Consider the Diophantine equation
$$F(x,y)=0$$
in integers $x,y$ and where $F$ is a polynomial with integer
coefficients. We shall assume that $F$ is irreducible over $\bbbc.$
The equation $F=0$ then represents a geometrically 
irreducible algebraic curve, which we denote by $C.$
Denote by $g$ the genus of $C$ and 
the number of branches at $\infty$ by $s$.
From a well-known theorem of Siegel [Si]
it follows that if $s+2g-2>0$, then the number of integer solutions 
to $F(x,y)=0$ is finite. Recently P.Corvaja and U.Zannier [CZ] gave a very surprising
alternative proof of this fact using W.M.Schmidt's subspace theorem.
Unfortunately, neither proof of Siegel's theorem gives an algorithm to actually
solve the general equation $F(x,y)=0$. Only in very special cases this
is possible. For example, A.Baker's method of linear forms in logarithms
allows one to solve equations of the form $y^q=f(x)$ for any given
$q$ and any $f\in\bbbz[x]$ having three or more distinct zeros
(see for example [ST, Ch 6]). Yu.Bilu [Bi]
studied necessary conditions for the applicability of Baker's method and
found several new instances where $F(x,y)=0$ can be solved in principle.

In this paper we take up an old paper of Runge [Ru, 1887] where equations of a 
particular kind are solved. As introduction consider the equation $F(x,y)=0$.
Let $d$ be the total degree of $F$ and
denote the sum of all terms of total degree $d$ in $F$ by $F_0$.
Suppose that $F_0$
factors as a product of two non-constant relatively prime factors $F_0=G_0H_0$.
This is called {\it Runge's condition}.
The branches at infinity of $C$ either correspond to $G_0$ or to $H_0$.
Runge's idea was to construct a polynomial $P(x,y)\in\bbbz[x,y]$, non-constant
on $C$, such that $P(x,y)\to0$
as we move to infinity along one of the branches corresponding to $G_0$.
For sufficiently large $x,y$ the integer points on these branches should
satisfy $P(x,y)=0$ because $P$ assumes integral values at these points.
We can then find them by elimination with $F(x,y)=0$.
Similarly we deal with the branches corresponding to $H_0$.
This idea is described in the introduction of Runge's paper, so
an algorithm to solve the equation is in principle there. In
the present paper we shall turn Runge's idea into an actual algorithm that can 
be fed to a computer. 

In [Sch] we find a generalisation of Runge's idea if one considers
weighted degrees. We shall present this generalisation in a slightly
different language using Newton polygons. 
Let us write $F=\sum_{m,n}f_{m,n}x^my^n$ where the summation extends over
a finite set of integer pairs. For each pair $m,n$ with $f_{m,n}\ne0$
we draw a 
rectangle with vertices $(0,0)$ and $(m,n)$ in the plane. The
{\it Newton-polygon} of $F$ is the convex hull of these rectangles. 
We denote it by $N_F$. The edges of $N_F$, not contained in the coordinate
axes, are called the {\it slopes} of $N_F$. To every slope $E$ of $N_F$
we can associate the sum $F_E(x,y)=\sum_{(m,n)\in E}f_{m,n}x^my^n$.
The theorem of Runge and Schinzel, presented in a slightly different form,
reads as follows.

\begin{theorem}[Runge, Schinzel]\label{runge}
Let the notation be as above. Suppose that the Newton polygon of $F$ has
either two distinct slopes, or one slope $E$ and $F_E$ factors into two
nonconstant, relatively prime polynomials in $\bbbz[x,y]$. Then the
equation $F(x,y)=0$ has finitely many solutions.
\end{theorem}

Using the proof of this theorem it is possible to give explicit upper
bounds for size of the solutions $(x,y)$. This is done in
[HS] or [W]. However, the upper bounds are so 
large, even for small parameters, that using them for an exhaustive search
on $x,y$ is impossible in practice.

It is the goal of the present paper to give a practical algorithm that
actually finds the solutions if the coefficients and degrees of $F$ are not
prohibitively large. We believe that a pleasant feature of our algorithm is,
that we do not use Puiseux series (only truncated power series)
and that we work entirely over $\bbbq$. 

For the sake of completeness we formulate Runge's Theorem in a
more algebraic geometric language and which also works in number
fields. We learnt this formulation from an informal note by
Yu.Bilu. Let $C$ be a smooth, connected algebraic curve defined over a
number field $K$. Let $S$ be a finite set of places of $K$, including the
infinite ones and let ${\cal O}_S$ be the ring of $S$-integers.
Fix a function $f\in C(K)$. A point $P\in C(K)$
will be called $S$-integral with respect to $f$ if $f(P)\in {\cal O}_S$.
The Galois group $Gal(\overline{K}/K)$ acts on the set of poles of $f$.
Denote by $\Sigma$ the set of orbits under this action. 

\begin{theorem}[Runge]. Assume that $|\Sigma|>|S|$. Then the $S$-integral
points of $C$ are effectively bounded.
\end{theorem}

For example, when $K=\bbbq$, $f$ is the $x$-coordinate and $S$ consists 
of the place at $\infty$, we see that we must have $|\Sigma|>1$, i.e.
the set of poles of $x$ must exist of at least two Galois orbits.
This is precisely
the factorisation condition discussed earlier. Instead of taking one function
$f$ we could have taken a finite set of functions $f_1,\ldots,f_t$. 
In our case over $\bbbq$ we would take $f_1=x,f_2=y$. 

Using ideas of Sprindzuk, Bombieri [Bo] found an interesting extension of
Runge's theorem. Let $s$ be a positive integer. A point $P\in C(K)$
is called $s$-integral if $|f(P)|_v>1$ for at most $s$ places $v$ of $K$.

\begin{theorem}[Bombieri, Sprindzuk]
Assume that $|\Sigma|>s$. Then the $s$-integral points of $C$ are effectively
bounded.
\end{theorem}

\section{Preparations}
To start with, we assume that the Newton
polygon of $F$ has a slope $E$ which is neither vertical nor
horizontal (we call this a {\it tilted slope}).
The remaining case, when the Newton polygon is a
rectangle, will be dealt with at the end of this section.

When there is only one slope, we assume that the associated
polynomial $F_E$ is a product of two relatively prime polynomials.
This is called the {\it Runge assumption}.
When there are only two slopes and one is horizontal and the other
$E$ we interchange $x,y$. We then get a new polynomial $F$
whose Newton polygon has a vertical slope and a tilted slope $E$.
After having made this change, if necessary, we are now in the
position that $F_E$ factors into two relatively prime polynomials
whose degrees in $y$ are strictly less than $\deg_y(F)$. 

Suppose that the points on the slope $E$ satisfy $ax+by=w$
where $a,b$ are relatively prime integers and $w$ is some integer.
We define the {\it weight} of a monomial $x^my^n$ by $am+bn$.
The weight of a polynomial $P$ is the maximum of the weights of
the monomials occurring in $P$. Notation: $w(P)$. In particular
we have that $w=w(F)$. 

Let us denote the factorisation of $F_E$ by $F_E=G_0H_0$.
We now give an algorithm to solve $F(x,y)=0$. 
First we consider the real points of $F(x,y)=0$. 
Let $x_2$ be the largest positive zero of the ${\rm discr}_y(F)$
and $x_1$ the smallest. Take $x_2=x_1$ to be specified later
if there are no real zeros. For any $x>x_2$ the equation
$F(x,y)=0$ has a fixed number, say $r$, of real solutions that we
denote by $y_1,y_2,\ldots,y_r$. We call the functions $y_i(x)$
the positive real branches of the curve $F=0$. Note that for any
branch $y_i(x)$ there is a
real zero $\alpha_i$ of $F_0(1,\alpha)$
such that $y_i(x)/x^{b/a}\to\alpha_i$ as $x\to\infty$. 
As a consequence of the Runge assumption we can find a
polynomial $P_i(x,y)$ with integer coefficients, not divisible
by $F$, such that $P_i(x,y_i(x))\to0$ as $x\to\infty$. The
construction of $P_i$ will be carried out in the next section.
Here we conclude the algorithm. Choose a positive parameter 
$\tau$. Let $x^+(i)$ be the largest positive zero of the resultant
of $F(x,y)$ and $P_i(x,y)+\tau$ with respect to $y$. Let $x^-(i)$
be the largest positive zero of the resultant of $F$ and $P_i(x,y)-\tau$.
Let $X_i$ be the maximum of $x_2,x^+(i),x^-(i)$ where we ignore $x_2$
if it has not been defined yet.
Then we know that for all $x>X_i$ we have 
$|P_i(x,y_i(x))|<\tau$. Suppose we have an integer point $(x,y)$ on
the $i$-th branch with $x>X_i$. Then $P_i(x,y)$ assumes an integral
value $a$ with $|a|<\tau$. We find all such $(x,y)$ simply by 
solving the simultaneous systems $F(x,y)=0,P_i(x,y)=a$ for
all integers $a$ with $|a|<\tau$.  

We carry out these steps for each positive real branch of $F=0$.
After that we have found all integer solutions $(x,y)$ with $x>\max_i X_i$.
Next we should consider the case $x<x_1$. For any such $x$
the equation $F(x,y)=0$ has a fixed number $r'$ of real solutions
$y_i(x)$ which we call the negative branches. For each such $i$ we
construct a function $P_i(x,y)$ and proceed as above.

In practice the number of distinct $P_i$ is smaller than the
number of actual real branches because one polynomial may vanish
on several asymptotic branches of the curve.

Finally, we promised to give an algorithm in the case when
the Newton polygon has no tilted slopes, i.e. it is a rectangle.
In that case $F$ contains a term $ax^my^n$ where $m=\deg_x(F)$
and $n=\deg_y(F)$. When $F(x,y)=0$ with $x$ large, the value of $y$ 
will be close to a zero of $F_2(y)=\lim_{x\to\infty}x^{-m}F(x,y)$. 
We now simply use the above algorithm by taking $F_2(y)$ as 
polynomial whose value tends to zero as we let $(x,y)$ follow
a branch of $F=0$ with $x\to\infty$.

\section{Construction of the vanishing functions}
Let notation be as in the introduction. Suppose we want to
compute the integral points on one of the positive real branches
$y/x^{b/a}\to \alpha$ of $F=0$. In this section we construct
a polynomial whose value on this branch tends to 0 as we let $x\to\infty$. 

Perform the
following change of variables, $x\to 1/t^b,y\to
\eta/t^a$. We obtain
$$F(x,y)={1\over t^{w(F)}}f(t,\eta)={1\over t^{w(F)}}(f_0+tf_1+t^2f_2+
\cdots)$$
where $f_i\in\bbbz[\eta]$ and $f_0=G_0(1,\eta)H_0(1,\eta)$. 
We now perform a Hensel lift of $f_0(\eta)=G_0(1,\eta)H_0(1,\eta)$
to a factorisation in $(\bbbq[[t]])[y]$ of the form
$$f(t,\eta)=(g_0+g_1t+\cdots)(h_0+h_1t+\cdots)$$
where $g_0=G_0(1,\eta),h_0=H_0(1,\eta)$, the degrees of $g_i,\ i>0$
are strictly less than the degree of $g_0$ and the degrees of
$h_i,\ i>0$ are strictly less than the degree of $h_0$. 

Notice that $\tilde{g}:=g_0(y)+g_1(y)t+\cdots=0$
is an analytic curve which,
for small $t$, contains a subset of the branches of $F=0$. 
Suppose our particular branch is among this union of branches.
We now construct a polynomial $P$ which vanishes on the branches
of $\tilde{g}=0$ as $x\to\infty$. Consider a polynomial 
$P\in\bbbz[x,y]$ with unknown coefficients and such that
$\deg_y P<\deg_y(F)$ and $w(P)\le N$ for some integer $N$ to be
specified. We rewrite
$P(1/t^a,\eta/t^b)=t^{-w(P)}p(t,\eta)$.
We choose our coefficients such that $p\mod{\tilde{g}}
=O(t^{N+1})$. 
The number $N$ is chosen in such a way that the number of coefficients
of $P$ exceeds the number of equations following from the 
constraint. 

\begin{lemma} The vector space
$$\{P\in\bbbq[x,y]\ |\ \deg_y P<\deg_y F,\ w(P)\le N\}$$
has $\bbbq$-dimension at least 
$N^2/(2\delta_x\delta_y)+N/(2\delta_x)+N/(2\delta_y)+C_1(F)$ if 
$N<w(F)$ and dimension at least $N\deg_y(F)/a+C_2(F)$ if $N\ge w(F)$.
Here $C_1,C_2$ depend only on $F$, not on $N$. 
\end{lemma}

As to the number of equations provided by
$p\mod{\tilde{g}}=O(t^{N+1})$, a priori we have $N\deg_{\eta}
(g_0)$ conditions. However, this number would in general be too
large. Fortunately we have the following additional consideration.
Let $\zeta$ be a primitive $a$-th root of unity. Replacing
$t$ by $\zeta t$ and $\eta$ by $\eta\zeta^b$ does not
change $P(1/t^a,\eta/t^b)$. Hence $p(t,\eta)$
changes by a factor $\zeta^{-b}$. A similar remark holds for
$\tilde{g}$. Hence, after computation of $p\mod{\tilde{g}}$,
the only terms that occur transform with the same character under
our substitution. So the actual number of constraints is
at most $N\deg_{\eta}(g_0)/a+1$.

Because $\deg_{\eta}(g_0)<\deg_y(F)$ we have that $N\deg_y(F)/a+C_2$
exceeds $N\deg_{\eta}(g_0)/a$ for sufficiently large $N$. Hence
there exists a polynomial $P$ of the required type.

\section{Example 1}
Consider the equation
$$F(x,y):=y^6 - 2y^5 - 4y^2x^4 + 17yx^2 + 4x - 18=0.$$
The highest degree part is given by $y^6-4y^2x^4$. We now 
replace $x\to 1/t,y\to \eta/t$ to get
$$f(t,\eta):=4t^5 - 18t^6 + 17t^3\eta - 4\eta^2 - 2t\eta^5 + \eta^6=0.$$
We shall be interested in Hensel lifts of the factorisation
$$\eta^6-4\eta^2=(\eta^2-2)(\eta^2+2)\eta^2$$
up to order $4$. We get
$$f(t,\eta)=g_1g_2g_3$$
where
\begin{eqnarray*}
g_1&=&\eta^2-2-\eta t-t^2/2+15\eta t^3/8+O(t^4)\\
g_2&=&\eta^2+2-\eta t-t^2/2+19\eta t^3/8+O(t^4)\\
g_3&=&\eta^2-17\eta t^3/4+O(t^4).
\end{eqnarray*}
First we determine a polynomial that vanishes on the branches
given by $g_1=0$. Let
$$P(x,y)=(ax+by)(y^2-2x^2)+py^2+qxy+rx^2+kx+ly+m$$
where $a,b,p,q,r,k,l,m$ are numbers to be determined. 
Define
$$p(t,\eta)=t^3P(1/t,\eta/t).$$
Then,
\begin{eqnarray*}
p\mod{g_1}&=&(2b + 2p + r + a\eta + q\eta)t\\
&& +(a/2 + k + 3b\eta/2 + l\eta + p\eta)t^2\\
&& +(-13b/4 + m + p/2 -15a\eta/8)t^3+O(t^4).
\end{eqnarray*}
This remainder vanishes up to order $4$ if and only if
$$\matrix{2b+2p+r=0 & a+q=0\cr
a/2+k=0 & 3b/2+l+p=0\cr
-13b/4+m+p/2=0 & -15a/8=0.\cr}$$
One solution is given by
$$a=0,b=2,p=-3,q=0,r=2,k=0,l=0,m=8.$$
So we get the desired function
$$P_1=2y(y^2-2x^2)+2x^2-3y^2+8.$$
The integer points lying on $g_2=0$ do not correspond to any
real branches extending to infinity because $y^2+2x^2$
has no real factors. Finally,
the function $y^2$ vanishes on the branches given by $g_3=0$
as $t\to0$. 

Let us now solve the equation $F(x,y)=0$. First
of all the real roots of the discriminant of $F$ with respect
to $y$ lie in between $1$ and $1.25$. 
First we deal with the branches given by $g_1=0$.
The real zeros of the resultant of $F$ and $P_1+1$ with
respect to $y$ lie between $-4.1$ and $3.2$.
The real zeros of the resultant of $F$ and $P_1-1$ are between
$-3.8$ and $3$. Hence, for the solutions on the branches given
by $g_1=0$ we have that $|P_1(x,y)|<1$ whenever $x\ge4$ or $x\le-5$.
In other words, we have $P_1(x,y)=0$ for such points. The $y$-resultant
of $F$ and $P_1$ has no integer zeros $x$.

As we said we can ignore the factor $g_2$ because it does not correspond
to ant branches. 
Finally we consider the branches given by $g_3=0$.
The resultant of $y^2+1$ and $F$ has no real zeros, the 
resultant of $y^2-1$ and $F$ has its real zeros between $1$ and $2$.
So $y^2<1$ for all integer points on the branches given by
$g_3=0$. In other words, $y^2=0$ for such points. Since $y=0$
and $F(x,y)=0$ imply $4x-18=0$, there are no integer solutions.

We are left to check the remaining values of $x$ between $4$ and $-5$.
So we check $F(k,y)=0$ for integer solution for $k=-4,-3,\ldots,2,3$.
It turns out that there are no integer solutions.

\section{Example 2}
Consider the equation
$$F(x,y):=y^4 + 2y^3 - 9x^2y^2 + 2xy - 15x - 7=0.$$
The highest degree part is given by $y^4-9x^2y^2=(y-3x)(y+3x)y^2.$
First we define
$$f(t,\eta)=t^4F(1/t,\eta/t).$$
Then factor $f$ up to order $4$ as
$$f(t,\eta)=g_1g_2g_3$$
where
\begin{eqnarray*}
g_1&=&\eta-3+t-t^2/18-13t^3/54+O(t^4)\\
g_2&=&\eta+3+t+5t^2/18-13t^3/54+O(t^4)\\
g_3&=&\eta^2-2\eta t^2/9+5t^3/3+O(t^4).
\end{eqnarray*}
Clearly $P_1:=y-3x+1$ vanishes on the branch given by $g_1=0$
and $P_2:=y+3x+1$ vanishes on the branch given by $g_2=0$.
A straightforward computation as in the previous section
shows that $P_3:=2y^3+15y^2$ vanishes on the branches given
by $g_3=0.$

We summarise the resultant computations in the following table.
\begin{center}
\begin{tabular}{|l|c|c|}
\hline
Resultant$_y$ & $x_{\rm min}$ & $x_{\rm max}$\\
\hline
$F,F_y$ & -1.03 & 0.78\\
$F,P_1+1$ & none & none \\
$F,P_1-1$ & -1.46 & 1.52\\
$F,P_2+1$ & -1.15 & 1.10\\
$F,P_2-1$ & -0.86 & 0.75\\
$F,P_3+1$ & -2.18 & 2.12\\
$F,P_3-3$ & -7.83 & 2.10\\
\hline
\end{tabular}
\end{center}
For the branches given by $g_1=0$ and $g_2=0$ we see that
$|P_1(x,y)|<1$ and $|P_2(x,y)|<1$ for all integer points
with $|x|\ge2$. In other words, $P_1(x,y)=P_2(x,y)=0.$
The system $P_1=F=0$ does not give integer solutions,
neither does $P_2=F=0.$

For the branches given by $g_3=0$ we see that $-1<P_3(x,y)<3$
for all $x$ with $x\le -8$ and $x\ge3.$ Hence $P_3(x,y)=0,1,2.$
Combining these possibilities with $F=0$ again does not yield
any integral solutions. 

It remains to check all points with $-7\le x\le 2.$ When $x=-1$
we get the solutions $y=-4,-1,1,2.$

\section{Example 3}
Consider the equation
$$F(x,y):=(y^2 - x^3)(y^2 - 2x^3) + 2x^5 - 9xy - 3=0.$$
This is an example which satisfies the Runge condition with respect
to the weight $2m+3n$ for a monomial $x^my^n.$
Introduce
$$f(t,\eta)=t^{12}F(1/t^2,\eta/t^3).$$
We can Hensel lift this to a factorisation
$$f(t,\eta)=g_1g_2$$
where
\begin{eqnarray*}
g_1&=&\eta^2-1 - 2t^2 - 4t^4 - 16t^6+O(t^7)\\
g_2&=&\eta^2-2 + 2t^2 + 4t^4 + 16t^6+O(t^7).
\end{eqnarray*}
We look for a function of the form
$$P:=ay^2 + by + cyx + px^3 + qx^2 + rx + s$$
which vanishes on the branches given by $g_1=0.$
Define
$$p(t,\eta):=t^6P(1/t^2,\eta/t^3).$$
Then 
\begin{eqnarray*}
p\mod{g_1}&=&(a+p)+c\eta t+(2a+q)t^2+b\eta t^3\\
&&+(4a+r)t^4+(16a+s)t^6+O(t^7).
\end{eqnarray*}
Note, by the way, that only terms $t^{2k},\eta t^{2l+1}$ occur.
To have this remainder vanish up to order $7$ we can take
$$a=-1,b=0,c=0,p=1,q=2,r=4,s=16$$
and we get the function
$$P_1=x^3-y^2+2x^2+4x+16.$$
Similarly we have
\begin{eqnarray*}
p\mod{g_2}&=&2a + p + c\eta t + (-2a+q)t^2+b\eta t^3\\
&&((-4a+r)t^4 + (-16a+s) t^6.
\end{eqnarray*}
To have this remainder vanish up to order $7$ we can take
$$a=1,b=0,c=0,p=-2,q=2,r=4,s=16$$
and we get the function
$$P_2=y^2-2x^3+2x^2+4x+16.$$
We now make our table of ranges of real zeros for
the various resultants.
\begin{center}
\begin{tabular}{|l|c|c|}
\hline
Resultant$_y$ & $x_{\rm min}$ & $x_{\rm max}$\\
\hline
$F,F_y$ & -1.50 & 8.23\\
$F,P_1+5$ & 15.0 & 31.7\\
$F,P_1-0.5$ & none & none \\
$F,P_2+5$ & 13.8 & 36.1\\
$F,P_2-1$ & none & none\\
\hline
\end{tabular}
\end{center}
So, when $x\le-2$ or $x\ge37$ we have either $-5< P_1(x,y)<0.5$
or $-5< P_2(x,y)<1$. First we solve $P_1(x,y)+k=0,F(x,y)=0$
for $k=0,1,2,3,4.$ There are no integer solutions. Then we
solve $P_2(x,y)+k=0,F(x,y)=0$ for $k=0,1,2,3,4$. Again there
are no solutions. Finally we solve $F(k,y)=0$ for $-1\le k\le 36$.
We find the solution $x=2,y=3.$

\section{References}
\begin{itemize}
\item[[Bi]] Yu.Bilu, Effective analysis of integral points on algebraic
curves, Israel J. Math. 90 (1995), 235-252.
\item[[Bo]] E.Bombieri, On Weil' s `Th\'eor\`eme de d\'ecomposition', 
Am.J.Math 105 (1983), 295-308.
\item[[CZ]] P.Corvaja, U.Zannier, A subspace theorem approach to
integral points on curves. C.R. Math. Acad. Sci. Paris 334(2002), 267-271.
\item[[HS]] D.L.Hilliker, E.G.Strauss, Determination of bounds for
the solutions to those diophantine equations that satisfy the 
hypotheses of Runge's theorem, Trans.Amer. Math.Soc. 280 (1983), 637-657.
\item[[Ru]] C.Runge, \"Uber ganzzalige L\"osungen von Gleichungen
zwischen zwei Ver\"anderlichen, J.reine angew. Math 100(1887), 425-435.
\item[[Sch]] A.Schinzel, An improvement of Runge's theorem on diophantine
equations, Comment.Pontificia Acad.Sci, 20 (1969), 1-9.
\item[[Si]] C.L.Siegel, \"Uber einige Anwendungen Diophantischer Approximationen,
Abh. Preuss. Akad. Wiss. Phys.-Math. Kl 1 (1929).
\item[[ST]] T.N.Shorey, R.Tijdeman, {\it Exponentaial Diophantine Equations},
Cambridge 1986.
\item[[W]] P.G.Walsh, A quantitative version of Runge's theorem on
Diophantine equations, Acta Arith. 62(1992), 157-172.
\end{itemize}

\end{document}